\documentclass[a4paper,11pt]{amsart}

\newcommand{\col}[1]{\textcolor{black}{#1}}

\usepackage{epsfig}

\usepackage{amsmath,amssymb,textcomp,graphicx}
\usepackage{mathrsfs}
\usepackage{times}
\usepackage[T1]{fontenc}
\usepackage[english]{babel}
\usepackage{rawfonts}
\usepackage{color}
\usepackage[latin1]{inputenc}
\def\cprime{$'$}

\newcommand{\R}{\mathbb{R}}

\newcommand{\N}{\mathbb{N}}

\newcommand{\C}{\mathbb{C}}

\renewcommand{\P}{\mathcal{P}}

\newcommand{\hA}{\hat{\mathcal{A}}}
\newcommand{\Id}{\mathrm{Id}}
\newtheorem{defi}{Definition}

\newtheorem{theorem}[defi]{Theorem}
\newtheorem{lemma}[defi]{Lemma}
\newtheorem{prop}[defi]{Proposition}
\newtheorem{cor}[defi]{Corollary}

\title[Squared eigenfunctions are generically linearly independent]{The squares of the Laplacian-Dirichlet eigenfunctions are generically linearly independent}
\author{Yannick Privat}
\author{Mario Sigalotti}
\address{Institut \'Elie Cartan de Nancy, UMR 7502
Nancy-Universit\'e - CNRS - INRIA, B.P. 239 , Vand\oe
uvre-l\`es-Nancy Cedex France}
\date{\today}

\begin{document}

\begin{abstract}
The paper deals with the genericity of domain-dependent spectral properties of the Laplacian-Dirichlet operator. In particular we prove that, generically, the squares of the eigenfunctions form a free family. We also show that the spectrum is generically non-resonant. 
The results are obtained by applying global perturbations of the domains 
and exploiting analytic perturbation properties.
The work is motivated by two applications: an existence result for 
the problem of maximizing the rate of exponential decay of a damped membrane and an approximate controllability result for the bilinear Schr\"odinger equation. 
\end{abstract}



\maketitle
\section*{Introduction}
Genericity is a measure of how much robust and frequent a property is. It  enjoys, therefore, a deep-rooted success in control theory, where 
a generic behavior 
is, roughly speaking, the expected behavior 
of systems 
involving
physical quantities whose value can only be approximated.

A paradigmatic example of generic properties 
in control theory is the controllability of a finite-dimensional linear system 
\begin{equation}\label{linear}
\dot x=Ax+Bu,\ \  x\in \R^n,\  \  u\in \R^m.
\end{equation}
It is well known, and the proof simply follows from the Kalman criterion,
 that for every choice of the positive integers $n$ and $m$ a generic linear system of type (\ref{linear}) is controllable. More precisely, the set of pairs $(A,B)$
 for which (\ref{linear}) is controllable
is open and dense in the product of the spaces of $n\times n$ and 
$n\times m$ matrices. (See, for instance, \cite{sontag-book}.)

When a control system involves partial differential equations,
conditions guaranteeing its controllability, observability or stabilizability 
can often be stated in terms of the eigenvalues or eigenspaces of some
linear operator (typically, the leading term of the evolution operator). 
In this paper we are mainly interested in conditions depending on the domain on which the control system of partial differential equations 
is defined. 
The genericity of some relevant conditions for control applications
has already been considered and proved in the general field of partial differential equations
(e.g., the simplicity of the eigenvalues of the Laplacian-Dirichlet
operator proved in \cite{micheletti-Perturbazione,ule} and applied 
in the control framework in
 \cite{zuazua-lions-hydroelastic}). Others, due to their technical nature,  
need to be tackled by specific arguments. 
This has led to the development of several tools 
for studying the  genericity with respect to the domain of control-related properties of partial differential operators. 
Without seeking exhaustiveness, 
let us mention the works by 
Lions and Zuazua \cite{lions-zuazua1996} and 
Ortega and Zuazua \cite{ortega-zuazua-stokes}
on the Stokes system,
those by Ortega and Zuazua \cite{ortega-zuazua-plate,ortega-zuazua-plate-addendum} on the plate equation, the paper by 
Chitour, Coron and Garavello \cite{CCG} on the heat and wave equations 
and the recent work by Beauchard, Chitour, Kateb and Long \cite{BCKL} on the 
Schr\"odinger equation.

The scope of this paper is to prove the  genericity with respect to the domain of some properties of the Laplacian-Dirichlet operator issuing  from 
control theory and optimization among which, in particular, the linear independence of the squared eigenfunctions. 

In doing so we propose 
a technique that, we believe, has a wider range of applicability, going beyond the conditions studied here and adaptable to different operators. 
The difference between our approach and those usually 
adopted 
is that we focus less on local infinitesimal variations of the domain and more 
on global, long-range perturbations. 
In order to get genericity results from this kind of perturbations 
 we have to rely on analytic-dependence properties 
for the eigenvalues and eigenfunctions of 
the Laplacian-Dirichlet operator with respect to analytic perturbations of the domain. (It should be stressed, however, that analytic perturbation theory applies to a much larger range of operators.)
The idea of proving genericity through global perturbations is clearly not new, being intrinsically contained in analytic perturbation theory.
Our work has actually been inspired 
by a paper by Hillairet and Judge \cite{hillairet-judge},
where the authors prove, using global perturbations, the generic simplicity of the eigenvalues 
of the Laplacian-Dirichlet operator on planar polygons with at least four vertices. 
The argument in \cite{hillairet-judge}, however, relies on 
the existence, in the class of interest, of domains having 
 simple spectrum.
The difficulty of extending the proof of \cite{hillairet-judge} to 
show the generic linear independence of the squared eigenfunctions on smooth domains is
that examples of 
smooth domains having the desired property are not handily available. 
One kind of domain on which the property can be easily checked is given by orthotopes.
However, many results on spectral stability
when non-smooth domains are
approximated  by smooth ones are known (see, in particular, 
the works by Arendt and Daners
\cite{Arendt-Daners} \col{and Bucur \cite{bucur}} where uniform stability  of the eigenfunctions is studied) and imply the existence, for every $n\in\N$, of a smooth domain $\mathcal{R}_n$ whose first $n$ eigenfunctions have linearly independent squares. 
In order to propagate by global analytic perturbation the property satisfied by $\mathcal{R}_n$ one can use, for instance, exponential flows
of vector fields (even a narrow family of vector fields is enough to 
generate a full orbit of domains, see \cite{agrachev-caponigro}). 
One has, however,  to take care of the possible crossing 
of the analytically depending eigenvalues. 
In order to do so, one should select analytic paths 
along which the first $n$ eigenvalues are simple. This problem is related to the Arnol\cprime d conjecture (see \cite{arnold-modes-and-quasi-modes,verdiere-commentari})
and has been solved  
by  Teytel in \cite{Teytel}. 
Teytel's result, recalled in Proposition~\ref{teytelB}, is crucial for 
the proposed perturbation technique (Theorems~\ref{MainTheo} and \ref{MainTheo-topology}).

Let us conclude this introduction by describing the motivating applications of the properties that we consider.
The generic linear independence of the squared eigenfunctions
has been conjectured in dimension two 
by H\'ebrard and Henrot in  \cite{Hebrard-Henrot-Optimal}, 
where the authors consider 
 the problem of stabilizing with the 
 largest possible decay rate
 a membrane fixed at its boundary using  a damping acting on a portion of the membrane of fixed area.
The existence and uniqueness of the solution for this problem can be deduced from the linear independence of the squared 
eigenfunctions of the Laplacian-Dirichlet operator 
on the domain filled by the membrane. 
(See Section~\ref{applications} for more details.)

It should be noticed that whether such linear independence is not only generic but rather always true is still an open question. 
A negative result by Mahar and Willner \cite{Mahar-Willner}
on the squared eigenfunctions of a Sturm-Liouville operator
justify a cautious stance toward a conjecture saying that the linear independence should always hold true. 

Linear independence of the squared eigenfunctions 
appears quite naturally also in the study of the controllability of the bilinear Schr\"odinger equation. In this context, indeed,
non-resonance conditions on the spectrum of the uncontrolled  Schr\"odinger 
operator are often required (see, for instance, \cite{CMSB}).
Since the $k^\mathrm{th}$ eigenvalue $\lambda_k^\varepsilon$ of
$-\Delta+\varepsilon V:H^2(\Omega)\cap H^1_0(\Omega)\to L^2(\Omega)$ is analytic with respect to $\varepsilon$ and satisfies
$$\left. \frac d {d\varepsilon}\right|_{\varepsilon=0}\lambda_k^\varepsilon=\int_\Omega V(x)\phi_k(x)^2 dx$$
where $(\phi_n)_{n\in\N}$ is a complete system of eigenfunctions of $-\Delta$ (see \cite{Albert}), then the linear independence 
of the family $(\phi_n^2)_{n\in\N}$ clearly plays a role in the study of 
the size of the family of potentials $V$ for which the spectrum has some prescribed property. 

Another application discussed in Section~\ref{applications}  
corresponds to the case where the uncontrolled  Schr\"odinger 
operator is defined by a potential well, i.e., $V=0$ and $\Omega$ is free.  We show in this case that, generically with respect to $\Omega$, 
no nontrivial linear combination with rational coefficients of the eigenvalues of $-\Delta$ annihilates. We deduce from this fact and the results in 
\cite{CMSB} a generic approximate controllability property for the Schr\"odinger equation.

Properties about the non-annihilation of linear combinations of eigenvalues play a role also in other domains. Let us mention, for instance,  the recent work by Zuazua on switching systems in infinite dimension \cite{zuazua-switching}, where the condition that the sums of two different pairs of eigenvalues of the Laplacian-Dirichlet operator are different is used to prove null-controllability of the heat equation using switching controls.

The paper is organized as follows: in Section~\ref{theor} we introduce some definitions and notations and we prove the main abstract results of the paper (Theorems~\ref{MainTheo} and \ref{MainTheo-topology}). We conclude the section by deducing from the abstract results some specific generic conditions; in particular, we obtain the generic linear independence of the squared eigenfunctions of the Laplacian-Dirichlet operator. In Section~\ref{applications} we propose two applications of  these generic properties to the stabilization of vibrating membranes and to the  controllability of the Schr\"odinger equation.

\par {\bf Acknowledgments.} We would like to thank 
Yacine Chitour, Antoine Henrot, Pier Domenico Lamberti  
and Enrique Zuazua 
for several fruitful discussions and advices.

\section{Generic properties by global perturbations}\label{theor}
\subsection{Notations and abstract genericity result}\label{notations}
Throughout the paper, $d$ denotes an integer larger than or equal to two and $\N$ the set of positive integer numbers, while $\N_0=\{0\}\cup \N$.

\par Given a Lipschitz domain $\Omega\subset \R ^d$, we denote by $(\lambda_n^\Omega)_{n\in \N}$ the nondecreasing sequence of eigenvalues of the Laplacian-Dirichlet operator 
$$-\Delta:H^2(\Omega)\cap H^1_0(\Omega)\to L^2(\Omega)$$ 
counted according to their multiplicity. As it is well known, it is always possible to choose an orthonormal basis of $L^2(\Omega)$ made of eigenfunctions of the Laplacian-Dirichlet operator. In the sequel any such choice will be denoted by $(\phi_n^\Omega)_{n\in \N}$ with  $\phi_n^\Omega$ corresponding to the eigenvalue 
$\lambda_n^\Omega$. 
We will identify $\phi_n^\Omega$ with its extension to zero outside $\Omega$.
\par We define the class of domains $\Sigma_m$ as the set of open \col{connected} subsets 
of $\R ^d$  with $\mathcal{C}^m$ boundary. 
By $D_m$ we denote the subset of $\Sigma_m$ 
of {\it $\mathcal{C}^m$ topological balls}, i.e., those open subsets  $\Omega$ of $\R ^d$ such that there exists 
a $\mathcal{C}^m$-diffeomorphism of $\R^d$ transforming the unit ball in $\Omega$. 
Similarly, 
we define $D_{0,1}$ as the orbit of the unit ball by bi-Lipschitz homeomorphisms of $\R^d$. 

It is well known that $\Sigma_m$ and $D_m$, endowed with the $\mathcal{C}^m$ topology
inherited from that of $\mathcal{C}^m$-diffeomorphisms,  are complete metric spaces (see \cite{Micheletti}). In particular, they are Baire spaces.
\par Let us recall that, given a Baire space $X$, a residual set (i.e. the intersection of countably many open and dense subsets) is dense in $X$. A boolean function $\P:X\to \{0,1\}$ is said to be \textit{generic} in $X$ if there exists a residual set $Y$ such that every $x$ in $Y$ satisfies property $\P$, that is, $\P(x)=1$.
\par A sequence of open domains $(\Omega_n)_{n\in\N}$ is said to 
{\it compactly converge} to a domain $\Omega$ if for every compact set $K\subseteq \Omega \cup \overline{\Omega}^c$, there exists $n_K\in \N$ such that for all $n\geq n_K$, $K\subseteq \Omega_n \cup \overline{\Omega_n}^c$.
\par In the sequel of the paper, we make use several times of the following result, whose proof can be found in \cite[Theorem 7.3]{Arendt-Daners}.
\begin{prop}\label{stable}
Let $n\in \N$ and fix a Lipschitz domain $\Omega\subset \R^d$ such that $\lambda_1^\Omega,\dots,\lambda_n^\Omega$ are simple. 
Let $\Omega_k$ be a sequence of Lipschitz domains compactly converging to $\Omega$ and such that $\cup_{k\in \N}\Omega_k$ is bounded.
Then $\lambda _j^{\Omega_k}\to \lambda_j^\Omega$ and, therefore, $\lambda_j^{\Omega_k}$ is simple for every $j=1,\ldots,n$, for $k$ large enough. Moreover, up to a sign in the choice of $\phi_j^{\Omega_k}$, $\phi_j^{\Omega_k}\to 
\phi_j^{\Omega}$ in $L^\infty(\R^d)$, as $k$ goes to infinity, for $j=1,\ldots ,n$.
\end{prop}
\par Another result playing a crucial role in our argument is the following 
\col{proposition. (See \cite[Theorem 6.4]{Teytel}.)}
\begin{prop}\label{teytelB}
\col{Let $m > 2$ and $\Omega_0$, $\Omega_1$ be two  
domains 
in $\Sigma_m$
that are $\mathcal{C}^m$-differentiably isotopic.
Then 
there exists an analytic curve $[0,1]\ni t\mapsto Q_t$ of $C^m$-diffeomorphisms such that $Q_0$ is equal to the identity, 
$Q_1(\Omega_0)=\Omega_1$ and  
every domain $\Omega_t=Q_t(\Omega_0)$
 has simple spectrum for $t$ in the open interval $(0,1)$. }
\end{prop}
Teytel 
deduces the proposition stated above \col{in the case where $\Omega_0$ and $\Omega_1$ are $\mathcal{C}^m$-differentiably isotopic to the unit $d$-dimensional ball} from a more general result, namely \cite[Theorem B]{Teytel}. 
\col{His argument applies also, without modifications, to pairs of domains belonging to the same isotopy class.} 
\col{Theorem B in \cite{Teytel}}   guarantees the existence of an analytic path of simple-spectrum operators among any elements of a family of operators satisfying 
a {\it strong Arnold hypothesis} on their eigenvectors 
(see also \cite{arnold-modes-and-quasi-modes,verdiere-commentari}). For this reason we expect that our method could be adapted to other situations.

We are ready to prove the following theorem 
on generic properties among topological balls. 
\begin{theorem}\label{MainTheo}
Let $F_{n}:\R ^{n(n+1)}\longrightarrow \R$, $n\in \N$, be a sequence of analytic functions.
For every $n\in \N$, 
we say that a Lipschitz domain $\Omega$ satisfies property $\P_n$ if 
$\lambda_1^\Omega, \dots,\lambda_n^\Omega$ are simple and if 
there exist $n$ points $x_1,\ldots,x_{n}$ in $\Omega$ and a choice 
$\phi_1^\Omega,\ldots,\phi_n^\Omega$ of the first $n$ eigenfunctions 
of the Laplacian-Dirichlet operator on $\Omega$ such that
\begin{equation}\label{Fn}
F_n (\phi _1^\Omega (x_1),\ldots , \phi _{n}^\Omega(x_1),\ldots,\phi_1^\Omega (x_n),\ldots,\phi_n^\Omega(x_n),\lambda _1^\Omega,\ldots,\lambda_n^\Omega )\neq 0.
\end{equation}
Assume that, for every $n\in \N$, there exists $\mathcal{R}_n\in D_{0,1}$ satisfying property $\P_n$. Then, for every $m\in \N\cup \{+\infty\}$, a generic $\Omega\in D_m$ satisfies $\P_n$ for every $n\in \N$.
\end{theorem}
\begin{proof}
Fix $m\in \N\cup \{+\infty\}$. 
Define, for every $n\in\N$,  the set of domains
$$
\mathcal{A}_n=\{\Omega\in D_m\mid 
\Omega\mbox{ satisfies }\P_n\}.
$$
We shall fix $n\in\N$ and prove that each $\mathcal{A}_n$ is open and dense in $D_m$.
\par Let us first prove that $\mathcal{A}_n$ is open. Fix $\Omega\in \mathcal{A}_n$, a choice of eigenfunctions $\phi_1^\Omega,\ldots,\phi_n^\Omega$ and 
$n$ points $x_1,\ldots,x_n\in \Omega$ such that (\ref{Fn}) holds true. Suppose by contradiction that there exists a sequence $(\Omega_k)_{k\in\N}$ in $D_m\setminus \mathcal{A}_n$ that converges to $\Omega$. 
Notice that the convergence in $D_m$ implies compact convergence in the sense recalled above.
Proposition~\ref{stable} thus implies that, for a choice of $\phi_j^{\Omega_k}$, $j=1,\ldots,n$, one has
\begin{eqnarray*}
\lefteqn{\lim_{k\to +\infty}F_n (\phi _1^{\Omega_k} (x_1),\ldots ,\phi_n^{\Omega_k}(x_n),\lambda _1^{\Omega_k},\ldots,\lambda_n^{\Omega_k})=} \\
& & F_n (\phi _1^{\Omega} (x_1),\ldots ,\phi_n^{\Omega}(x_n),\lambda _1^{\Omega},\ldots,\lambda_n^{\Omega})\neq 0.
\end{eqnarray*}
This contradicts the assumption that $\Omega_k\notin\mathcal{A}_n$ for every $k\in\N$.
\par 
We prove now the density of $\mathcal{A}_n$. Notice that, without loss of generality, $m>2$. Fix $\Omega\in D_m$.
Let $\mathcal{R}_n$ be as in the statement of the theorem, that is, 
$\mathcal{R}_n\in D_{0,1}$ and satisfies property $\P_n$. 
Notice that $\mathcal{R}_n$ can be approximated by a sequence of domains in $D_m$ in the sense of the compact convergence. Therefore, by applying the same argument as above, 
we deduce that
there exists $\widetilde{\mathcal{R}}_n\in D_m$ satisfying $\P_n$. Choose $\phi_j^{\widetilde{\mathcal{R}}_n}$, $j=1,\dots,n$, and $x_1,\dots,x_n\in \widetilde{\mathcal{R}}_n$ such that
$$
F_n (\phi _1^{\widetilde{\mathcal{R}}_n}(x_1),\ldots , \phi _{n}^{\widetilde{\mathcal{R}}_n}(x_n),\lambda _1^{\widetilde{\mathcal{R}}_n},\ldots,\lambda_n^{\widetilde{\mathcal{R}}_n} )\neq 0.
$$
 
We now apply Proposition~\ref{teytelB} with $\Omega_0=\widetilde{\mathcal{R}}_n$ and $\Omega_1=\Omega$. 
We deduce that, for $m>2$, 
there exists an analytic curve $[0,1]\ni t\mapsto Q_t$ of $C^m$-diffeomorphisms such that $Q_0$ is equal to the identity, 
$Q_1(\widetilde{\mathcal{R}}_n)=\Omega$ and  
every domain $\Omega_t=Q_t(\widetilde{\mathcal{R}}_n)\in D_m$
 has simple spectrum for $t$ in the open interval $(0,1)$. 
Due to standard analytic perturbation theory (see \cite{Kato}),
$\lambda_k^{\Omega_t}$ are analytic functions of $t$
 and there exists a choice of $\phi_j^{\Omega_t}$, $j=1,\dots,n$, $t\in[0,1]$, such 
that $\phi_j^{\Omega_t}\circ Q_t$ 
varies analytically with respect to $t$ in 
$\mathcal{C}^m(\widetilde{\mathcal{R}}_n)$. 
In particular,
$$
t\mapsto F_n (\phi _1^{\Omega_t}(Q_t(x_1)),\ldots , \phi _{n}^{\Omega_t}(Q_t(x_n)),\lambda _1^{\Omega_t},\ldots,\lambda_n^{\Omega_t} )
$$
is an analytic real-valued function. Since its value at $t=0$ is different from zero, then it annihilates only for finitely many $t\in[0,1]$. 

Hence, as required, $\Omega$ can be approximated arbitrarily well in $D_m$ by 
an element of $\mathcal{A}_n$. 
\end{proof}

Let us turn our attention to domains that are not necessarily
topological balls. 
\col{
The extension of Theorem~\ref{MainTheo} works 
along similar lines, once a deformation argument is used to transfer 
each property $\P_n$
from the set of topological balls to the desired isotopy class of domains. 
}

\begin{theorem}\label{MainTheo-topology}
\col{Let $(F_n)_{n\in\N}$, $(\P_n)_{n\in\N}$  and 
$(\mathcal{R}_n)_{n\in\N}$
be as in the statement of Theorem~\ref{MainTheo}. 
 Then, for every $m\in \N\cup \{+\infty\}$, a generic $\Omega\in \Sigma_m$ satisfies $\P_n$ for every $n\in \N$.
}
\end{theorem}
\begin{proof}
Fix $m\in \N\cup \{+\infty\}$. 
\col{Thanks to Theorem~\ref{MainTheo}, 
a generic $\hat\Omega\in D_m$ satisfies $\P_n$ for every $n\in \N$.}
Fix one such $\hat\Omega$ and notice that, in particular, the spectrum
$(\lambda_n^{\hat\Omega})_{n\in\N}$ is simple.

Define, for every \col{$n\in \N$},  the set 
$$
\col{\hA_n=\{\Omega\in \Sigma_m\mid 
\Omega\mbox{ satisfies }\P_n\}.}
$$

The openness of \col{$\hA_n$} in $\Sigma_m$ can be proved 
following exactly the same argument used in the proof of Theorem~\ref{MainTheo} to show that each $A_n$ is open in $D_m$. 

We are left to prove that \col{$\hA_n$} is dense in $\Sigma_m$. \col{Without loss of generality $m>2$.}
Take $\Omega\in \Sigma_m$. 
Let $B$ be an open ball of $\R^d$ containing $\Omega$. 
By eventually shrinking $B$, we can assume that
$\partial B\cap \partial \Omega$ contains at least one point $p$. 
Up to a change of coordinates, we can assume that $B$ is centered at the origin and $p=(0,\dots,0,1)$. 

Consider a smooth vector field on $\R^d$ satisfying
$$V(x_1,\dots,x_d)=\left\{\begin{array}{ll}
\left( \begin{array}{c}x_1x_d\\ \vdots\\ x_{d-1}x_d\\ x_d^2-\frac{x_1^2+\cdots+x_d^2+1}2 \end{array}\right)&\mbox{if $x_1^2+\cdots+x_d^2< \rho$}\\
0&\mbox{if $x_1^2+\cdots+x_d^2>\rho+1$}
\end{array}\right.
$$
for some $\rho>1$. The behavior of $V$ in a neighborhood of the unit ball is
represented in Figure~\ref{figV}. Notice that $V$ is complete, since it vanishes outside a compact set. 
 
\begin{figure}[h!]
\begin{center}
\input{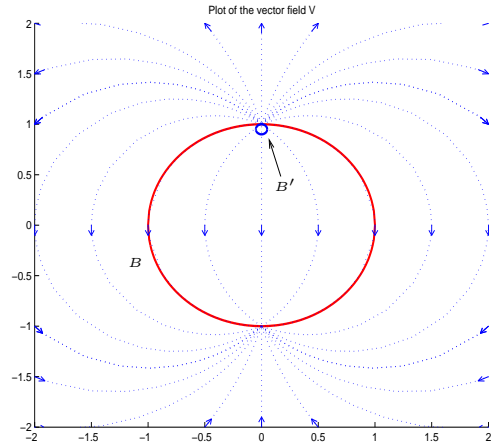}
\caption{\label{figV}Phase portrait of the vector field $V$. 
}
\end{center}
\end{figure}

By construction $V$  is everywhere tangent to $\partial B$. 
The ball $B$ is therefore invariant for the flow of $V$. 
Notice that the points $p$ and $-p$ are the only zeros of $V$ in $\overline{B}$ and that $x_d$ is strictly decreasing along all trajectories of $V$ staying in $B$.
Therefore, $p$ is a repulsive equilibrium for $V$ restricted to $\overline{B}$ and $-p$ an attractive one.  

Notice that, since $\Omega\in\Sigma_m$, then its boundary has finitely many components and therefore there exists  a ball $B'$ contained in $\Omega$  
such that $p\in \partial B'$. 
Notice, moreover, that the differential of $V$ at $\pm p$ is $\pm \Id$.  
Then, for every $x\in B$, $e^{-t V}(x)$ belongs to $B'$ for every $t$ larger than some $t_x\in\R$. 
We deduce that $e^{t V}(B')$ 
compactly converges to $B$
as $t$ tends to infinity. Since $B'\subset \Omega$,  then $e^{t V}(\Omega)$ 
compactly converges to $B$ as well as $t\to+\infty$.

Consider an analytic path $t\mapsto Q_t$  of $\mathcal{C}^m$-diffeomorphisms of $\R^d$ such that $Q_0=\Id$ and 
$Q_1(B)=\hat\Omega$, whose existence can be deduced from 
Proposition~\ref{teytelB}. Then 
$$\Omega_t=Q_{\frac{2\arctan t}\pi}\circ e^{t V}(\Omega)\mbox{ compactly converges to }\hat\Omega\mbox{ as $t\to+\infty$}.$$
Moreover, \col{each $\Omega_t$ is isotopic to $\Omega$.} 
\col{It follows from Proposition~\ref{stable}
that we can fix  $t$ large enough in such a way that $\Omega_t$ verifies $\P_n$. 
Proposition~\ref{teytelB} implies that there exists an analytic path of domains $s\mapsto \tilde \Omega_s$ such that $\tilde\Omega_0=\Omega$, $\tilde\Omega_1=\Omega_t$ and the spectrum of the Laplacian-Dirichlet
operator on $\tilde \Omega_s$ is simple for every $s\in(0,1)$. 
}

\col{Hence, as in the proof of Theorem~\ref{MainTheo}, we can 
deduce that $\tilde \Omega_s$ satisfies $\P_n$ for all but finitely many $s\in[0,1]$. In particular, $\Omega$ is in the closure of $\hA_n$.}
\end{proof}

\par 
\subsection{Consequences of the abstract results}
In this section, we present two corollaries of Theorem~\ref{MainTheo-topology} 
showing that (i) the squares of the Laplacian-Dirichlet eigenfunctions are generically
linearly independent and (ii) the Laplacian-Dirichlet spectrum is generically non-resonant.

Recall that a finite or infinite sequence of real numbers 
is said to be {\it non-resonant} if every nontrivial rational linear combination of 
finitely many of  
its elements 
is different from zero.

\par 
In order to verify that the squares of the Laplacian-Dirichlet
eigenfunctions on a suitably chosen $d$-orthotope   are linearly independent,
we prove the following
technical result.
\begin{lemma}\label{lemmaIndpdt}
Let $\varphi$ belong to $\mathcal{C}^\infty([0,+\infty),\R)$, $N$ be a positive integer and $(a_1,\dots,a_N)$
be a sequence of pairwise distinct positive real numbers. 
Assume that there exist $l_0\in \N_0$ and $l_1\in \N$ such that $\varphi ^{(l_0+p l_1)}(0)\neq 0$ for every $p=0,\dots,N-1$. Then, the functions 
$\varphi (a_1\cdot),\dots,\varphi (a_N\cdot)$ are linearly independent on every right-neighborhood of zero. 
\end{lemma}
\begin{proof}
We are interested in finding all the $N$-tuples
$(\gamma_1,\ldots,\gamma_N)\in\R^N$ such that $\sum_{k=1}^N\gamma_k\varphi
(a_k\cdot )=0$ in a right-neighborhood of zero. Differentiating this relation 
$l_0+p l_1$ times yields the relation
$\sum_{k=1}^N\gamma_k a_k^{l_0+p l_1}\varphi ^{(l_0+p l_1)}(a_k \cdot)=0$. Evaluating such relation at zero 
for $p=0,\dots,N-1$, 
we obtain a system of $N$ linear equations in the $N$ variables $\gamma_1,\dots,\gamma_N$. 
Since $(a_j^{l_1(i-1)})_{1\leq
i,j\leq N}$ is a Vandermonde matrix, the determinant $\delta_N$ of 
the $N\times N$ matrix underlying such a
system writes
\begin{eqnarray*}
\delta_N & = & \det \left((a_j^{l_1(i-1)})_{1\leq i,j\leq
N}\right)\prod_{k=1}^N a_k^{l_0}\varphi^{(l_0+k l_1)}(0)\\
 & = & \prod_{1\leq i<j\leq
N}(a_j^{l_0}-a_i^{l_0})\prod_{k=1}^Na_k^{l_0}\varphi^{(l_0+k l_1)}(0)\neq 0.
\end{eqnarray*}
This concludes the proof of the lemma.
\end{proof}

\begin{prop}\label{propRectangle}
Let $(\mu_1,\ldots,\mu_d)$ be a non-resonant sequence of positive real numbers
 and $\mathcal{R}$ be the
$d$-orthotope $\prod_{i=1}^d(0,\mu_i\pi)$. 
Then, the Laplacian-Dirichlet
eigenvalues of $\mathcal{R}$ are simple 
and the squares of the Laplacian-Dirichlet eigenfunctions are linearly
independent.
\end{prop}
\begin{proof}
Let us prove the lemma by induction on $d\geq 1$. 

If $d=1$, then $\mu_1$ is any positive real number and 
the squared eigenfunctions of the Laplacian-Dirichlet operator on 
$\mathcal{R}$ are $(\sin^2(k \cdot/\mu_1))_{k\in\N}$.
The proposition follows then from Lemma~\ref{lemmaIndpdt}, 
applied to $\varphi(x)=\sin^2(x)$, $l_0=1$, $l_1=2$, and $a_k=k/\mu_1$. 

Let now $d$ be larger than one.
For $K=(k_1,\ldots,k_d)\in {\N}^d$, we write $K'=(k_1,\ldots,k_{d-1})$, $\bar K=k_d$ and we denote by 
$f_{K}$ the (un-normalized) Laplacian-Dirichlet eigenfunction 
$$
f_{K}(x_1,\ldots,x_d)=
\prod_{i=1}^d\sin \left(\frac{k_i x_i}{\mu_i}\right).
$$
Clearly, $f_K(x_1,\ldots,x_d)=f_{K'}(x_1,\ldots,x_{d-1}) \sin (\bar K x_d/\mu_d)$. 
Fix $I\subset \N^d$ finite and $\{\gamma_K\mid K\in I\}\subset \R$ such that
$\sum_{K\in I} \gamma_K f_{K}^2  \equiv 0$ on $\mathcal{R}$. 
Let $\bar I=\{\bar K\mid K\in I\}$. 
Then for every $(x_1,\dots,x_{d-1})\in \prod_{i=1}^{d-1}(0,\mu_i\pi)$
and every $x_d\in (0,\mu_d \pi)$ we have
$$\sum_{k\in\bar I}\left(\sum_{K\in I,\bar K=k}  \gamma_K f_{K'}(x_1,\dots,x_{d-1})^2\right) \sin^2 \left(\frac{k x_d}{\mu_d}\right)    = 0.$$
Therefore, applying again Lemma~\ref{lemmaIndpdt} to  
$\varphi(x)=\sin^2(x)$, we deduce that, for every $k\in\bar I$,
$$\sum_{K\in I,\bar K=k}  \gamma_K f_{K'}^2\equiv  0\ \ \ \ \mbox{on } \prod_{i=1}^{d-1}(0,\mu_i\pi).$$
The induction hypothesis implies that 
$\gamma_K$ is equal to zero for every $K$ such that $\bar K=k$. 
Since $k$ is arbitrary in $\bar I$, the proposition is proved. 
\end{proof}
We can now state the first corollary of Theorem~\ref{MainTheo-topology}.
\begin{cor}\label{SquareEigenRectangle}
Let $m\in \N \cup \{\infty\}$. Generically with respect to $\Omega\in \Sigma_m$, the squares of the Laplacian-Dirichlet eigenfunctions are linearly independent when restricted to any  measurable subset of $\Omega$ of positive measure.\end{cor}
\begin{proof}
First notice that $n$ functions $\varphi_1,\dots,\varphi_n$ defined on a domain $\Omega$ are linearly independent if and only if there exist 
$n$ points $x_1,\dots,x_n$ in  $\Omega$ such that
$$
 \det \left(\begin{array}{ccc}
\varphi_1(x_1) & \ldots &\varphi_n(x_1) \\
\vdots & & \vdots \\
\varphi_1(x_n) & \ldots & \varphi_n(x_n)
\end{array}\right)\ne 0.
$$

Apply Theorem~\ref{MainTheo-topology} with
$$
F_n(y_1,\ldots,y_{n(n+1)}) = \det \left(\begin{array}{ccc}
y_1 & \ldots & y_n \\
\vdots & & \vdots \\
y_{n^2-n+1} & \ldots & y_{n^2}
\end{array}\right),
$$
for $(y_1,\ldots,y_{n(n+1)})\in \R ^{n(n+1)}$, and \col{${\mathcal{R}}_n=\mathcal{R}$ for every $n\in \N$}, where 
$\mathcal{R}$ is the $d$-orthotope introduced in the statement of 
Proposition~\ref{propRectangle}.

Then for a generic $\Omega\in \Sigma_m$ the squares of the Laplacian-Dirichlet eigenfunctions are linearly independent on $\Omega$.
 Assume that there exists a measurable subset $\mathcal{O}\subset \Omega$ 
 of positive measure and $K$ constants $\gamma_1,\dots,\gamma_K$ 
 such that 
$\sum_{k=1}^K\gamma_k\phi_k^\Omega(x)^2=0$ on $\mathcal{O}$.
Recall now that  the hypo-analyticity of the
Laplacian operator implies that each eigenfunction is analytic
inside $\Omega$. 
 Hence $\gamma_1=\cdots=\gamma_K=0$. 
\end{proof}

Corollary~\ref{SquareEigenRectangle}
can be used to get generic spectral properties as in \cite[Section 6.3]{zuazua-switching}.

Another consequence of Theorem~\ref{MainTheo-topology} is the following corollary.

\begin{cor}\label{non-res}
Fix $m\in \N\cup\{\infty\}$, $k\in\N$ and $q=(q_1,\dots,q_k)\in\R^k\setminus\{0\}$.
Then, for a generic $\Omega\in \Sigma_m$ one has
\begin{equation}\label{q}
\sum_{l=1}^kq_l \lambda_l^\Omega \ne 0.
\end{equation}
In particular, a generic $\Omega\in \Sigma_m$ has non-resonant spectrum.
\end{cor}
\begin{proof}
Let $\mathcal{R}$ be a  
 $d$-orthotope defined as in the statement of 
Proposition~\ref{propRectangle}.  

 We denote by $\Gamma$ the subset of $\partial \mathcal{R}$ 
 defined by 
$$
\Gamma=\{(x_1,\ldots,x_d)\in \partial \mathcal{R}\mid x_d=\mu_d\pi\}.
$$
Consider a perturbation $\mathcal{R}^t:=(\Id+tV)(\mathcal{R})$ of the domain $\mathcal{R}$, with $t$ small and $V$ a smooth vector field whose 
support is compact and does not intersect $\partial \mathcal{R}\setminus \Gamma$. Then, it is well known (see, e.g. \cite{Henrot-Pierre,Simon,Soko}) that, since the Laplacian-Dirichlet eigenvalues of $\mathcal{R}$ are simple, the shape derivative of $\lambda_l^{\mathcal{R}} $ along $V$ is
defined as
$$
\langle d\lambda_l^{\mathcal{R}} ,V\rangle =\left.\frac{d}{dt}\lambda_l^{\mathcal{R}^t}\right|_{t=0}=
-\int_{\Gamma}\left(\frac{\partial \phi_l^{\mathcal{R}}}{\partial \nu}\right)^2(V\cdot \nu) d\sigma ,
$$
where $\nu$ denotes the outward normal to $\mathcal{R}$ and $d\sigma$ the $(d-1)$-dimensional surface element. 
By hypothesis $\nu=(0,\dots,0,1)$ on $\Gamma$, so that $V\cdot \nu$ is equal to $v_d$, the $d^\mathrm{th}$ component of $V$.  
Notice, moreover, that
$$\frac{\partial \phi_l^{\mathcal{R}}}{\partial \nu}=c_l f_{K_l}$$
for some nonzero constant $c_l\in \R$ (defined up to sign) and some $K_l\in \N^{d-1}$, where $f_{K_l}$ is defined as in the proof of Proposition~\ref{propRectangle}.

\par Let $q=(q_1,\dots,q_k)\in\R^k\setminus\{0\}$ 
\col{and introduce $G:\Omega \mapsto \sum_{l=1}^kq_l \lambda_{l}^\Omega $}. Differentiating \col{$G$} at $\Omega =\mathcal{R}$ along a vector field $V$ chosen as above yields, 
$$
\col{\langle dG ,V\rangle =-\int_{\Gamma}\sum_{l=1}^k q_l c_{l}^2 f_{K_{l}}^2(x_1,\dots,x_{d-1}) v_d\, d\sigma.}
$$
Due to Proposition~\ref{propRectangle}, 
$$\col{\sum_{l=1}^k q_l c_{l}^2 f_{K_{l}}^2|_\Gamma}$$ 
is not everywhere zero on $\Gamma$.  
Thus, it is possible to choose \col{$V$} for which \col{$
\langle dG ,V\rangle \ne 0$.}

The conclusion follows by applying  Theorem~\ref{MainTheo-topology} 
with $F_n=1$ for $n\ne k$ and 
$F_k(y_1,\dots,y_{k(k+1)})=\sum_{i=1}^k q_i y_{k^2+i}$
and by taking  \col{${\mathcal{R}}_k=(\Id+t V)(\mathcal{R})$} for $t$ small enough. 
\end{proof}

\section{Applications to shape optimization and control theory}
\label{applications}

\subsection{Stabilization of a damped membrane}
We consider here a stabilization problem in $\R ^2$ and we are interested in proving the existence and uniqueness of solutions for a related shape optimization problem. 
More precisely, let us denote by $\Omega\subset \R ^2$ a 
domain belonging to $D_m$, $m\in\N\cup\{\infty\}$. 
Assume that the Laplacian-Dirichlet eigenvalues of $\Omega$ are simple.
\par We consider the problem of stabilizing a membrane fixed at the boundary $\partial \Omega$, thanks to a damping acting only on a subdomain $\omega$. 
Denote by $\chi_\omega$ the characteristic function of $\omega$.
The displacement $v$ of the membrane, in presence of a viscous damping of the type $2k\chi_{\omega}$, $k>0$, satisfies
\begin{equation}
\left\{\begin{array}{ll}
\frac{\partial ^2v}{\partial t^2}-\Delta v+2k\chi_\omega (x)\frac{\partial v}{\partial t}=0 & (t,x)\in (0,+\infty)\times \Omega \\
v(t,x)=0 & x\in \partial \Omega , \ t>0\\
v(0,x)=v_0(x) & x\in \Omega \\
\frac{\partial v}{\partial t}(0,x)=v_1(x) & x\in\Omega,
\end{array}
\right.
\end{equation}
where $v_0\in H^1_0(\Omega)$ and $v_1\in L^2(\Omega)$. 
This system is known to be exponentially stable if $\omega$ has positive measure and it is possible to define its exponential decay rate (which does not depend on the initial data). A natural question consists in looking for the largest decay rate once the area of $\omega$ is fixed. 
Such optimization problem is already quite difficult
in the
one-dimensional case  (see e.g. \cite{Cox-Zuazua}). 
For this reason H\'ebrard and Henrot in  \cite{Hebrard-Henrot-Optimal} 
introduce a simplified version of it by considering, instead of the decay rate, the quantity
\begin{equation}
J_N(\omega):=\inf_{1\leq n\leq N}\int_{\Omega}\chi_{\omega}(x)(\phi_n^\Omega(x))^2 dx,
\end{equation}
where $N$ is a given positive integer and $\phi_n^\Omega$ denotes, as in the previous sections, the $n^\mathrm{th}$ normalized Laplacian-Dirichlet eigenfunction.
\par Then, we are driven to study the following shape optimization problem
\begin{equation}\label{ShapeOptimPb}
\left\{
\begin{array}{l}
\min J_N(\omega)\\
\omega \in \mathcal{L}_\ell ,
\end{array}
\right.
\end{equation}
where $\mathcal{L}_\ell$ denotes the set of 
measurable subsets of $\Omega$ of measure $l$. 
It is convenient to identify subdomains of $\Omega$ with their characteristic functions, so that $\mathcal{L}_\ell$ is identified with
$$
\left\{a\in L^\infty(\Omega)\mid a(x)=0\mbox{ or }1\mbox{ a.e. and }\int_\Omega a(x)dx=\ell\right\}.
$$
The one-dimensional problem is completely solved in \cite{Hebrard-Henrot-Spillover}. In the same paper it is noticed that
the proof of existence and uniqueness of the optimum for 
(\ref{ShapeOptimPb}) can be easily adapted to the two-dimensional case 
under the generic hypothesis that the square of the  Laplacian-Dirichlet  eigenfunctions $\phi_1^\Omega,\dots,\phi_N^\Omega$ are linearly independent (see Corollary~\ref{SquareEigenRectangle}). 
Indeed, first the authors prove the existence of an optimum $a^*$ in a relaxed class. In order to prove that such a maximum is a characteristic function, they 
study the optimality conditions satisfied by $a^*$, by considering perturbations of $a^*$ with support in ${A}_\varepsilon :=\{x\in \Omega \mid \varepsilon \leq a^*(x)\leq 1-\varepsilon\}$, with a small $\varepsilon >0$. 
They can prove in this way the existence of 
$N$ real numbers $\alpha _1,\dots,\alpha_N$ such that $\alpha_1^2+\cdots+\alpha_N^2\ne 0$ and
$$
\sum_{k=1}^N\alpha _k \phi_k^\Omega(x)^2=\mbox{constant, for almost every } x\in {A}_\varepsilon.
$$
Then, because of the analyticity of the eigenfunctions and of the linear independence of their squares, ${A}_\varepsilon$ must have measure zero.

\begin{theorem}
Generically with respect to $\Omega\in D_m$, the optimization problem (\ref{ShapeOptimPb}) has a unique solution $\omega_N^*$.
\end{theorem}

\subsection{Controlled Schr\"odinger equation}
We apply in this section Corollary~\ref{non-res} in order to prove the generic approximate controllability of a bilinear Schr\"odinger equation of the type
\begin{equation}\label{sch_EQ}
\left\{\begin{array}{ll}
i\frac{\partial \psi}{\partial t}(t,x) = (-\Delta +u(t)W(x))\psi(t,x), & (t,x)\in  (0,+\infty)\times \Omega  \\
\psi(t,x)=0 & x\in \partial \Omega , \ t>0\\
\psi(0,x)=\psi_0(x) & x\in \Omega, 
\end{array}
\right.
\end{equation}
where $\Omega$ belongs to $\Sigma_m$ for some $m\in \N\cup\{\infty\}$, $W\in L^\infty(\Omega,\R)$,  
the control $u$ belongs to $L^\infty([0,+\infty),U)$
for some fixed measurable subset $U$ of $\R$ with nonempty interior, and
 $\psi_0\in L^2(\Omega,\C)$. 
 System (\ref{sch_EQ}) admits always a  mild solution
  $\psi\in \mathcal{C}([0,+\infty),L^2(\Omega,\C))$ in the sense of \cite{BMS}.

The control system (\ref{sch_EQ}) is said to be approximately controllable if 
for every $\psi_0,\psi_1\in L^2(\Omega,\C)$ and every $\varepsilon>0$ there exist a control $u\in L^\infty([0,+\infty),U)$ and a positive time $T$ such that 
the solution $\psi$ of (\ref{sch_EQ}) satisfies $\|\psi(T,\cdot)-\psi_1\|_{L^2(\Omega)}<\varepsilon$. 

It has been proved in \cite{CMSB} that (\ref{sch_EQ}) is approximately controllable if  the Laplacian-Dirichlet operator on $\Omega$ has non-resonant spectrum and 
\begin{equation}\label{couple}
\int_\Omega W(x)\phi_k^\Omega(x)\phi_{k+1}^\Omega(x)\,dx\ne 0\ \ \ \mbox{ for every $k\in \N$}.
\end{equation}
Corollary~\ref{non-res} ensures that the Laplacian-Dirichlet spectrum is generically non-resonant. On the other hand, the unique continuation property implies that, 
 for every $k\in\N$, the product $\phi_k^\Omega\phi_{k+1}^\Omega$ is a nonzero function on $\Omega$.  Therefore, for every $\Omega$ with non-resonant spectrum,
$\{W\in L^\infty(\Omega)\mid \mbox{(\ref{couple}) holds true}\}$ is 
residual in $L^\infty(\Omega)$.
Moreover, due to the continuity of the eigenfunctions stated in Proposition~\ref{stable}, for every $k\in \N$ the map
$$(\Omega,W)\mapsto \int_\Omega W(x)\phi_k^\Omega(x)\phi_{k+1}^\Omega(x)\,dx$$ 
is continuous with respect to the product topology of 
$\Sigma_m\times L^\infty(\R^d)$. 
As a consequence 
we obtain the following result.
\begin{prop}
Generically with respect to $(\Omega,W)\in \Sigma_m\times L^\infty(\R^d)$, endowed with the product topology,
 system (\ref{sch_EQ})
 is approximately controllable.
\end{prop}

\bibliographystyle{plain}
\bibliography{biblio}  

\end{document}